\documentclass[a4paper,11pt]{article}
\usepackage{latexsym}
\usepackage{amssymb}
\usepackage{amsfonts}
\usepackage{amsmath}
\usepackage{indentfirst}
\usepackage{graphicx}
\usepackage{epsfig}
\usepackage{tikz,ifthen}
\newtheorem{prop}{Proposition}[section]
\newtheorem{teor}{Theorem}[section]

\newtheorem{cor}{Corollary}[section]

\newcommand{\ninN}{n\in \mathbf{N}}

\newcommand{\cvd}{\hfill $\blacksquare$\bigskip}
\newcounter{indice}

\newcommand{\shape}[1]{
\setcounter{indice}{0};
\foreach \i in {#1} {
\addtocounter{indice}{1};
\foreach \x in {1,...,\i} {
\draw (\x-1,-\theindice+1) rectangle (\x,-\theindice);
}
}
}

\newcommand{\griddedshape}[1]{
\setcounter{indice}{0};
\foreach \i in {#1} {
\addtocounter{indice}{1};
\foreach \x in {1,...,\i} {
\draw[gray,dotted] (\x-1,-\theindice+1) rectangle (\x,-\theindice);
}
}
}

\newcommand{\skewshape}[1]{
\setcounter{indice}{0};
\foreach \i/\j in {#1} {
\addtocounter{indice}{1};
\foreach \x in {\j,...,\i} {
\draw (\x-1,-\theindice+1) rectangle (\x,-\theindice);
}
}
}

\newcommand{\permutation}[1]{
\setcounter{indice}{0};
\foreach \i in {#1}
\addtocounter{indice}{1};

\addtocounter{indice}{1}
\draw [help lines] (1,1) grid (\theindice,\theindice);

\setcounter{indice}{1};

\foreach \i in { #1 } {
\draw (\theindice+.5,\i+.5) [fill] circle (.2);
\addtocounter{indice}{1};
}
}

\date{}

\author{Mathilde \textsc{Bouvel}\footnote{Mathilde Bouvel conveys special acknowledgements to the \emph{Dipartimento di Sistemi e Informatica}
for the kind hospitality during her visit in May 2008 in which this research started.} \\ LaBRI UMR 5800, Universit\'e de Bordeaux and
CNRS,\\ 351, cours de la Lib\'eration, 33405 Talence cedex, France\\
\texttt{bouvel@labri.fr}\\
\and Luca \textsc{Ferrari} \\ Dipartimento di Sistemi e Informatica, Universit\`a degli Studi di Firenze, \\ Viale Morgagni 65, 50134 Firenze, Italy\\ \texttt{ferrari@dsi.unifi.it}}

\title{On the enumeration of $d$-minimal permutations}

\begin{document}

\maketitle

\begin{abstract}
We suggest an approach for the enumeration of minimal permutations
having $d$ descents which uses skew Young tableaux. We succeed in
finding a general expression for the number of such permutations
in terms of (several) sums of determinants. We then generalize the
class of skew Young tableaux under consideration; this allows in
particular to discover some presumably new results concerning
Eulerian numbers.
\end{abstract}

\section{Introduction}

This article deals with minimal permutations with $d$ descents
(also called $d$-minimal permutations here). This family of
permutations has been introduced in \cite{BoRo} in the study of
the whole genome duplication-random loss model of genome
rearrangement. In this context, genomes are represented by
permutations, and minimal permutations with $d = 2^p$ descents are
the basis of excluded patterns that describes the class of
permutations that can be obtained from the identity with cost at
most $p$.

In order to describe properties of this class of permutations, its
basis has been studied, and the first natural question to address
is to count how many excluded patterns it contains. In \cite{BP}
some partial results on the enumeration of minimal permutations
with $d$ descents have been obtained: namely, minimal permutations
with $d$ descents and of size $n$ have been enumerated by closed
formulas, for $n = d+1, d+2$ and $2d$ ($d+1$ and $2d$ being lower
and upper bounds for the size of a minimal permutation with $d$
descents -- see \cite{BP}). In \cite{MY}, further results on the
enumeration of minimal permutations with $d$ descents have been
obtained using multivariate generating functions, allowing in
particular to derive a closed formula enumerating those of size
$2d-1$ as well as some asymptotic results.

In this work we offer an alternative approach for the enumeration
of minimal permutations with $d$ descents, making extensive use of
a bijection between these permutations and a family of skew Young
tableaux. This gives a general formula for the number $p_{d+k,d}$
of minimal permutations with $d$ descents and of size $d+k$, as a
sum of determinants of matrices (Theorem \ref{principale}). This
expression for $k=d$ will specialize into a determinant expression
of Catalan numbers which is believed to be new. When specializing
it for $k=3$, it also allows us to give a closed formula for
$p_{d+3,d}$ (Theorem \ref{d+3}). Finally, the family of skew Young
tableaux under consideration has a natural generalization which is
investigated in Section \ref{section_generalization}.

\section{Preliminary definitions and results}\label{prel}

For any integer $n$, $S_n$ denotes the set of permutations of
$[1..n]$. A permutation $\sigma \in S_n$ will be represented
either by the word $\sigma(1) \ldots \sigma(n)$ or by the $n\times
n$ grid, where a cell contains a dot if and only if it is at
coordinates $(i, \sigma(i))$ for some $i \in [1..n]$.

The \emph{pattern involvement order} on permutations \cite{P} is
defined as follows. A permutation $\pi \in S_k$ is involved in (or
is a pattern of) $\sigma \in S_n$ when there exist integers $1\leq
i_1 < \ldots < i_k \leq n$ such that $\pi$ and $\sigma(i_1) \ldots
\sigma(i_k)$ are order-isomorphic sequences, \emph{i.e.} they are
such that $\pi(\ell) < \pi(m) \Leftrightarrow \sigma(i_{\ell}) <
\sigma(i_m)$ for all $\ell, m \in [1..k]$.

A \emph{descent} in a permutation $\sigma \in S_n$ is an integer $i \in [1..(n-1)]$ such that
$\sigma(i) > \sigma(i+1)$. Similarly, an \emph{ascent} is an integer $i \in [1..(n-1)]$ such that
$\sigma(i) < \sigma(i+1)$.

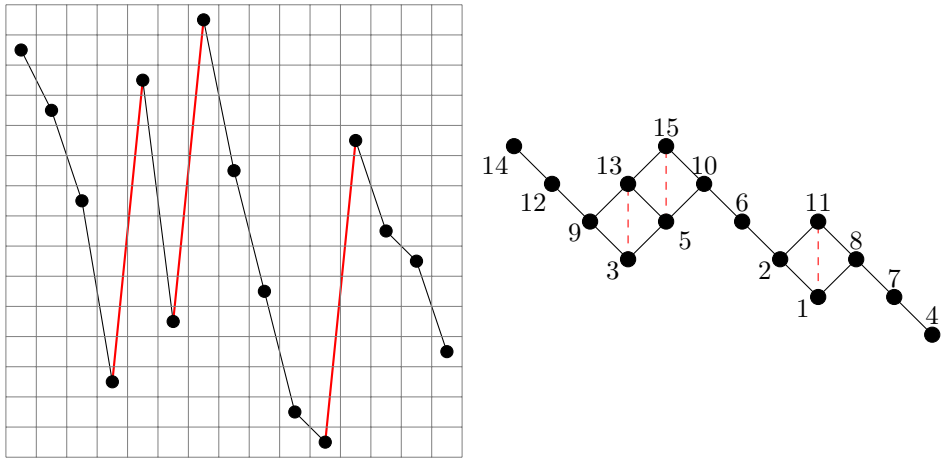
\begin{figure}[ht]
\begin{center}
\begin{tikzpicture}
\begin{scope}[scale=0.4]
\draw (1.5,14.5) -- (2.5,12.5);
\draw (2.5,12.5) -- (3.5,9.5);
\draw (3.5,9.5) -- (4.5,3.5);
\draw[thick, color=red] (4.5,3.5) -- (5.5,13.5);
\draw (5.5,13.5) -- (6.5,5.5);
\draw[thick, color=red] (6.5,5.5) -- (7.5,15.5);
\draw (7.5,15.5) -- (8.5,10.5);
\draw (8.5,10.5) -- (9.5,6.5);
\draw (9.5,6.5) -- (10.5,2.5);
\draw (10.5,2.5) -- (11.5,1.5);
\draw[thick, color=red] (11.5,1.5) -- (12.5,11.5);
\draw (12.5,11.5) -- (13.5,8.5);
\draw (13.5,8.5) -- (14.5,7.5);
\draw (14.5,7.5) -- (15.5,4.5);
\permutation{14,12,9,3,13,5,15,10,6,2,1,11,8,7,4}
\end{scope}
\end{tikzpicture}  \begin{tikzpicture}
\begin{scope}[scale=0.5]
\draw (1,0) node {~};
\draw[color=red, dashed] (4,5) -- (4,7);
\draw[color=red, dashed] (5,6) -- (5,8);
\draw[color=red, dashed] (9,4) -- (9,6);
\draw (1,8) [fill] circle (.2);
\draw (0.5,7.5) node {\small{$14$}};
\draw (2,7) [fill] circle (.2);
\draw (1.5,6.5) node {\small{$12$}};
\draw (3,6) [fill] circle (.2);
\draw (2.6,5.7) node {\small{$9$}};
\draw (4,5) [fill] circle (.2);
\draw (3.6,4.7) node {\small{$3$}};
\draw (4,7) [fill] circle (.2);
\draw (3.5,7.5) node {\small{$13$}};
\draw (5,6) [fill] circle (.2);
\draw (5.5,5.5) node {\small{$5$}};
\draw (5,8) [fill] circle (.2);
\draw (5,8.5) node {\small{$15$}};
\draw (6,7) [fill] circle (.2);
\draw (6,7.5) node {\small{$10$}};
\draw (7,6) [fill] circle (.2);
\draw (7,6.5) node {\small{$6$}};
\draw (8,5) [fill] circle (.2);
\draw (7.6,4.7) node {\small{$2$}};
\draw (9,4) [fill] circle (.2);
\draw (8.6,3.7) node {\small{$1$}};
\draw (9,6) [fill] circle (.2);
\draw (9,6.5) node {\small{$11$}};
\draw (10,5) [fill] circle (.2);
\draw (10,5.5) node {\small{$8$}};
\draw (11,4) [fill] circle (.2);
\draw (11,4.5) node {\small{$7$}};
\draw (12,3) [fill] circle (.2);
\draw (12,3.5) node {\small{$4$}};
\draw (1,8) -- (4,5);
\draw (4,7) -- (5,6);
\draw (5,8) -- (9,4);
\draw (9,6) -- (12,3);
\draw (3,6) -- (5,8);
\draw (4,5) -- (6,7);
\draw (8,5) -- (9,6);
\draw (9,4) -- (10,5);
\end{scope}
\end{tikzpicture}
\end{center}
\caption{The $11$-minimal permutation $\sigma = 14 \, 12\, 9\, 3\, 13\, 5\, 15\, 10\, 6\, 2\, 1\, 11\, 8\, 7\, 4$, and the corresponding poset. \label{fig:perm_and_poset}}
\end{figure}

A \emph{minimal permutation with $d$ descents}, or
\emph{$d$-minimal permutation}, of length $n$ is a permutation of
$S_n$ that is minimal in the sense of the pattern-involvement
relation for the property of having $d$ descents. In other words,
it is a permutation with $d$ descents such that, when removing any
of its entries and suitably renaming the remaining elements, the
resulting permutation of $S_{n-1}$ has $d-1$ descents. For
instance (see Figure \ref{fig:perm_and_poset}), the permutation
$\sigma = 14 \, 12\, 9\, 3\, 13\, 5\, 15\, 10\, 6\, 2\, 1\, 11\,
8\, 7\, 4$ is minimal with $11$ descents, as it has exactly $11$
descents and every permutation it involves as a pattern has at
most $10$ descents.

In \cite{BP}, minimal permutations with $d$ descents have been characterized as follows:

\begin{teor}
A permutation $\sigma$ is minimal with $d$ descents if and only if
it has exactly $d$ descents and its ascents $i$ satisfy the ``
diamond property", \emph{i.e.} are such that $2 \leq i \leq n-2$
and $\sigma(i-1) \sigma(i) \sigma(i+1) \sigma(i+2)$ forms an
occurrence of either the pattern $2143$ or the pattern $3142$.
\label{thmBP}
\end{teor}

As explained in \cite{BP}, this characterization allows to
represent $d$-minimal permutations by means of certain labelled
posets. These posets, labelled with the integers from $1$ to $n$,
are made of chains somehow linked by diamond-shaped structures
(corresponding to the ascents of the permutation). Figure
\ref{fig:perm_and_poset} shows an example of this one-to-one
correspondence. Notice that each of these labelled posets
represents a unique $d$-minimal permutation, whereas the
underlying unlabelled poset can be seen as representing a set of
$d$-minimal permutations (those that are in correspondence with a
legal labelling of the poset).

Posets and labelled posets in these families are in one-to-one
correspondence with skew Ferrers diagrams and skew Young tableaux
having special properties. These combinatorial objects have been
widely studied in the literature, in particular from an
enumerative point of view (see, for instance, the recent paper
\cite{BaRo}). In the following sections, we explicitly describe
the correspondence between unlabelled (resp. labelled) posets and
skew Ferrers diagrams (resp. skew Young tableaux), as well as some
enumerative results on these objects, and how they can be used for
our purposes.

\section{Connection with skew Young tableaux}

In order to explain how the poset representation of $d$-minimal permutations described in
the previous section can be conveniently interpreted by using skew
Young tableaux, we first need to recall some definitions.

\bigskip

An \emph{integer partition} of an integer $n$ is a sequence of positive integers
$\lambda = (\lambda_1, \ldots , \lambda_k)$ such that $\lambda_i
\geq \lambda_{i+1}$ for $1\leq i \leq k-1$ and $\sum_{i=1}^k
\lambda_i =n$. The integer $n$ is called the \emph{size} of the
integer partition, and we write $n = |\lambda|$. The number of
parts $k$ will be denoted by $k = \ell(\lambda)$ (this is also
called the \emph{length} of $\lambda$).

\begin{figure}[ht]
\begin{center}
\begin{tikzpicture}
\begin{scope}[scale=0.5]
\shape{8,6,3,3,2,1}
\end{scope}
\end{tikzpicture} \qquad \qquad
\begin{tikzpicture}
\begin{scope}[scale=0.5]
\shape{8,6,3,3,2,1}
\draw (0.5,-0.5) node {$23$};
\draw (1.5,-0.5) node {$22$};
\draw (2.5,-0.5) node {$21$};
\draw (3.5,-0.5) node {$19$};
\draw (4.5,-0.5) node {$16$};
\draw (5.5,-0.5) node {$12$};
\draw (6.5,-0.5) node {$9$};
\draw (7.5,-0.5) node {$5$};
\draw (0.5,-1.5) node {$20$};
\draw (1.5,-1.5) node {$17$};
\draw (2.5,-1.5) node {$13$};
\draw (3.5,-1.5) node {$8$};
\draw (4.5,-1.5) node {$4$};
\draw (5.5,-1.5) node {$1$};
\draw (0.5,-2.5) node {$18$};
\draw (1.5,-2.5) node {$14$};
\draw (2.5,-2.5) node {$6$};
\draw (0.5,-3.5) node {$15$};
\draw (1.5,-3.5) node {$10$};
\draw (2.5,-3.5) node {$3$};
\draw (0.5,-4.5) node {$11$};
\draw (1.5,-4.5) node {$2$};
\draw (0.5,-5.5) node {$7$};
\end{scope}
\end{tikzpicture}
\end{center}
\caption{The Ferrers diagram associated with the integer partition $\lambda = (8,6,3,3,2,1)$, and a Young tableau on this shape.\label{fig:Young_tableau}}
\end{figure}
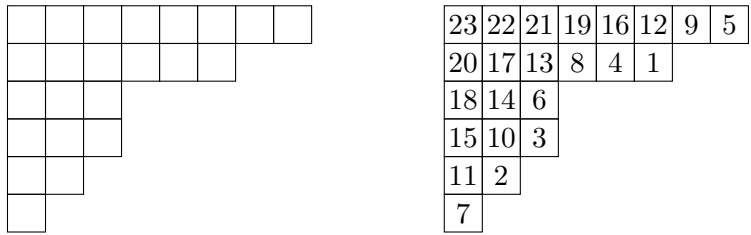

An integer partition $\lambda = (\lambda_1, \ldots , \lambda_k)$
can be represented by its \emph{Ferrers diagram}, which is
obtained by drawing $k$ rows of contiguous unit cells, from top to
bottom, such that row $i$ contains $\lambda_i$ cells, and with the
first cells of these $k$ rows vertically aligned. An example is
shown in Figure \ref{fig:Young_tableau}. We will also denote by
$\lambda$ the Ferrers diagram associated with the integer
partition $\lambda$. The size $|\lambda|$ obviously corresponds to
the number of cells of the Ferrers diagram, and the number of rows
is given by $\ell(\lambda)$.

For our purposes, a \emph{Young tableau} is a filling of a Ferrers
diagram $\lambda$ using distinct positive integers from $1$ to
$n=| \lambda |$, with the properties that the values are
(strictly) decreasing along each row and each column of the
Ferrers shape. This constitutes a slight departure from the
classical definition, which requires the word ``increasing"
instead of the word ``decreasing". However, it is clear that all
the properties and results on (classical) Young tableaux can be
translated into our setting by simply replacing the total order
``$\leq$" with the total order ``$\geq$" on $\mathbf{N}$. In
Figure \ref{fig:Young_tableau} a Young tableau of shape $\lambda =
(8,6,3,3,2,1)$ is shown. Like for Ferrers diagrams, the size of a
Young tableau is given by the number of its cells.

\bigskip

\bigskip

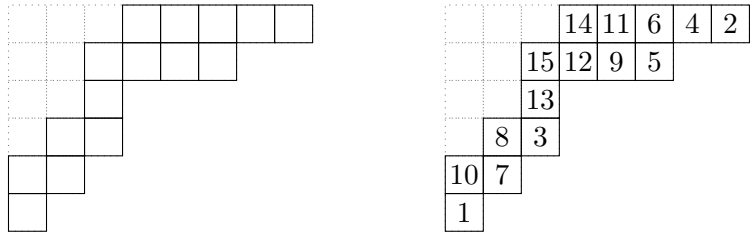
\begin{figure}[ht]
\begin{center}
\begin{tikzpicture}
\begin{scope}[scale=0.5]
\griddedshape{8,6,3,3,2,1}
\skewshape{8/4,6/3,3/3,3/2,2/1,1/1}
\end{scope}
\end{tikzpicture} \qquad \qquad
\begin{tikzpicture}
\begin{scope}[scale=0.5]
\griddedshape{8,6,3,3,2,1}
\skewshape{8/4,6/3,3/3,3/2,2/1,1/1}
\draw (3.5,-0.5) node {$14$};
\draw (4.5,-0.5) node {$11$};
\draw (5.5,-0.5) node {$6$};
\draw (6.5,-0.5) node {$4$};
\draw (7.5,-0.5) node {$2$};
\draw (2.5,-1.5) node {$15$};
\draw (3.5,-1.5) node {$12$};
\draw (4.5,-1.5) node {$9$};
\draw (5.5,-1.5) node {$5$};
\draw (2.5,-2.5) node {$13$};
\draw (1.5,-3.5) node {$8$};
\draw (2.5,-3.5) node {$3$};
\draw (0.5,-4.5) node {$10$};
\draw (1.5,-4.5) node {$7$};
\draw (0.5,-5.5) node {$1$};
\end{scope}
\end{tikzpicture}
\end{center}
\caption{The skew shape $(8,6,3,3,2,1) \setminus (3,2,2,1)$, and a skew Young tableau on this shape. \label{fig:skew_tableau}}
\end{figure}

The main definition we need in our work is that of a \emph{skew
Young tableau}. The definition can be given exactly as for a Young
tableau, with the only difference that the underlying shape
consists of a Ferrers diagram $\lambda$ with a Ferrers diagram
$\mu$ removed (starting from the top-left corner). Such a skew
shape is usually denoted $\lambda \setminus \mu$. We refer the reader to \cite{St1}
for the formal definition and some important facts concerning the
enumeration of skew Young tableaux. In Figure \ref{fig:skew_tableau} a skew Young
tableau of skew shape $(8,6,3,3,2,1) \setminus (3,2,2,1)$ is depicted. As before, the size of a skew Young tableau denotes its number of cells.

\bigskip

\bigskip

As announced at the beginning of the present section, we can
translate the poset representation of a $d$-minimal permutation
into a suitable skew Young tableaux.

\begin{prop}\label{bij} The set of $d$-minimal permutations of length $d+k$ is in
bijection with the set of skew Young tableaux whose skew shapes
$\lambda \setminus \mu$ satisfy $|\lambda \setminus \mu|=d+k$,
having $k$ rows and such that two consecutive rows have precisely
two columns in common.
\end{prop}

\emph{Proof.}\quad A $d$-minimal permutation of length $d+k$
consists of $k$ descending runs and, denoting with $a,b,c,d$ four
consecutive elements such that $a,b$ and $c,d$ belong to different
descending runs, then necessarily $a>b$, $c>d$, $a<c$ and $b<d$
(see Theorem \ref{thmBP}). Then, starting from a $d$-minimal
permutation $\pi$ of length $d+k$, one can construct a skew Young
tableau as follows: starting from the bottom, the $i$-th row of
the tableau consists of the elements of the $i$-th descending run
of $\pi$; moreover two consecutive rows are required to have
exactly two columns in common. The resulting tableau is skew Young
thanks to the above recalled diamond property of $d$-minimal
permutations.
\begin{flushright} $\blacksquare$ \end{flushright}

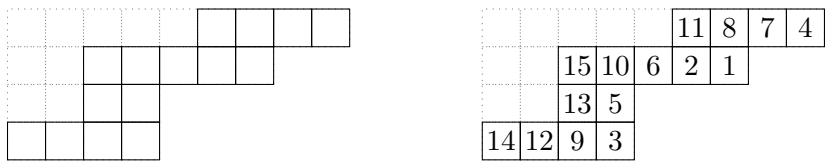
\begin{figure}[ht]
\begin{center}
\begin{tikzpicture}
\begin{scope}[scale=0.5]
\griddedshape{9,7,4,4}
\skewshape{9/6,7/3,4/3,4/1}
\end{scope}
\end{tikzpicture} \qquad \qquad \begin{tikzpicture}
\begin{scope}[scale=0.5]
\griddedshape{9,7,4,4}
\skewshape{9/6,7/3,4/3,4/1}
\draw (5.5,-0.5) node {$11$};
\draw (6.5,-0.5) node {$8$};
\draw (7.5,-0.5) node {$7$};
\draw (8.5,-0.5) node {$4$};
\draw (2.5,-1.5) node {$15$};
\draw (3.5,-1.5) node {$10$};
\draw (4.5,-1.5) node {$6$};
\draw (5.5,-1.5) node {$2$};
\draw (6.5,-1.5) node {$1$};
\draw (2.5,-2.5) node {$13$};
\draw (3.5,-2.5) node {$5$};
\draw (0.5,-3.5) node {$14$};
\draw (1.5,-3.5) node {$12$};
\draw (2.5,-3.5) node {$9$};
\draw (3.5,-3.5) node {$3$};
\end{scope}
\end{tikzpicture}
\end{center}
\caption{The skew shape corresponding to the underlying unlabelled poset of Figure \ref{fig:perm_and_poset}, and the skew Young tableau corresponding to the permutation $\sigma$ of Figure \ref{fig:perm_and_poset}. \label{fig:skew_shape_with2overlaps}}
\end{figure}

In Figure \ref{fig:skew_shape_with2overlaps} the skew Young tableau determined by the
permutation whose poset representation is given in Figure \ref{fig:perm_and_poset} is
shown.

\bigskip

\emph{Example.}\quad For $k=2$ and $d=4$, the set of $4$-minimal permutations of length $6$
is in bijection with the set consisting of all skew Young tableaux
of one of the following skew shapes:
\begin{itemize}
\item $(4,2) \setminus \emptyset$, \emph{i.e.} \begin{tikzpicture} \begin{scope}[scale=0.3] \shape{4,2} \end{scope} \end{tikzpicture},
\item $(4,3) \setminus (1)$, \emph{i.e.} \begin{tikzpicture} \begin{scope}[scale=0.3] \griddedshape{4,3} \skewshape{4/2,3/1} \end{scope} \end{tikzpicture},
\item $(4,4) \setminus (2)$, \emph{i.e.} \begin{tikzpicture} \begin{scope}[scale=0.3] \griddedshape{4,4} \skewshape{4/3,4/1} \end{scope} \end{tikzpicture}.
\end{itemize}

There are respectively $9$, $14$ and $9$ skew Young tableaux of these shapes, giving a total of $32$ $4$-minimal permutations of length $6$.

\bigskip

\emph{Remarks.}\quad Let $\pi$ be a $d$-minimal permutation and
suppose it has $k$ descending runs. Obviously this means that
$\pi$ has $k-1$ ascents and that $|\pi |=d+k$. Moreover, suppose
that $\lambda \setminus \mu$ is the skew shape associated with $\pi$.
Recall that $\ell(\lambda)$ and $\ell(\mu)$ denote the number of rows of $\lambda$ and $\mu$ respectively.
Finally, let $\ell_i$ be the length of the $i$-th descending run of
$\pi$. Some straightforward consequences of the above bijection
are the following:

\begin{enumerate}

\item The skew shape associated with $\pi$ is connected.

\item $\ell(\lambda )=k$.

\item Set $\lambda =(\lambda_1 ,\lambda_2 ,\ldots \lambda_k )$.
Then $\lambda_i =\sum_{j=1}^{k-i+1}\ell_j -2(k-i)$.

\item $\ell(\mu )<\ell(\lambda )$ and, more precisely, $\ell(\mu
)=(\ell(\lambda )-1)-\#(\textnormal{starting descending runs of
length 2 in $\pi$})$.

\item Set $\mu =(\mu_1 ,\mu_2 ,\ldots \mu_k )$. Then $\mu_i
=\lambda_{i+1}-2$.

\end{enumerate}

\section{Some enumerative results}

The main goal of the present section is to enumerate $d$-minimal
permutations of length $n$, with $d+1\leq n\leq 2d$. A general
result in this direction can be obtained by considering the above
described bijection with skew Young tableaux. In particular, an
interesting result due to Aitken is our starting point.

\begin{teor}\label{Ai}\textnormal{(\cite{A})} Let $\lambda \setminus \mu$ be a skew shape, with
$|\lambda \setminus \mu |=N$ and $\ell(\lambda ) = n$. Then,
the number $f^{\lambda \setminus \mu}$ of skew Young tableaux of shape
$\lambda \setminus \mu$ is
\begin{equation}\label{aitken}
f^{\lambda \setminus \mu}=N!\det\left( \frac{1}{(\lambda_i -\mu_j
-i+j)!}\right)_{i,j=1..n}.
\end{equation}
\end{teor}

Formula (\ref{aitken}) can be deduced from the well known
Jacobi-Trudi identity, as shown, for instance, in \cite{St2}.
Moreover, in such a formula the entry $(i,j)$ of the considered
matrix is intended to be 0 if the expression $\lambda_i -\mu_j
-i+j$ is negative.

The main result of this section is essentially a corollary of
Theorem \ref{Ai}, in the case in which the skew shape $\lambda \setminus \mu$ has the properties of Proposition \ref{bij}.

\begin{teor}\label{principale} Denote by $p_{d+k,d}$ the number of
$d$-minimal permutations of length $d+k$ (so that $1\leq k\leq
d$). Then
\begin{displaymath}
p_{d+k,d}=\sum_{a_1 ,a_2 ,\ldots ,a_k \geq 2\atop a_1 +a_2 +\cdots
+a_k =d+k}(d+k)!\cdot \det (A(a_1 ,\ldots ,a_k )),
\end{displaymath}
where $A(a_1 ,\ldots ,a_k )$ is the following matrix:
\begin{displaymath}
\left(
\begin{array}{ccccccc}
\frac{1}{a_1 !}&\frac{1}{(a_1 +a_2 -1)!}&\frac{1}{(a_1 +a_2 +a_3
-2)!}&\cdots&\cdots&\cdots&\frac{1}{(a_1 + \cdots +a_k -k+1 )!}
\\ 1&\frac{1}{a_2 !}&\frac{1}{(a_2 +a_3 -1)!}&\cdots&\cdots&\cdots&\frac{1}{(a_2 + \cdots +a_k -k+2 )!}
\\ \chi_{a_2 =2}&1&\frac{1}{a_3 !}&\ddots&&&\vdots
\\ 0&\chi_{a_3 =2}&1&\ddots&&&\vdots
\\ 0&0&\chi_{a_4 = 2}&\ddots&&&\vdots
\\ \vdots&\ddots&0&\ddots&&&\vdots
\\ \vdots&&&\ddots&&&\vdots
\\ \vdots&&&&&&\frac{1}{(a_{k-2}+a_{k-1}+a_k -2)!}
\\ \vdots&&\ddots&\ddots&\ddots&\ddots&\frac{1}{(a_{k-1}+a_k -1)!}
\\ 0&\cdots&\cdots&0&\chi_{a_{k-1}=2}&1&\frac{1}{a_k !}
\end{array}\right) .
\end{displaymath}
Here $\chi_P$ denotes the characteristic function of the property
$P$ (i.e., $\chi_P=1$ when $P$ is true and $\chi_P=0$ otherwise).
In other words, $A(a_1 ,\ldots ,a_k )$ is the $k\times k$ matrix
whose entries $a_{i,j}$ obey the following equalities:
\begin{eqnarray*}
a_{i,j}&=&\frac{1}{(a_i +\cdots +a_j +i-j)!},\qquad
\textnormal{when $i\leq j$,}
\\ a_{i,i-1}&=&1,
\\ a_{i,i-2}&=&\chi_{a_{i-1}=2},
\\ a_{i,j}&=&0,\qquad \textnormal{when $i>j+2$.}
\end{eqnarray*}
\end{teor}

\emph{Proof.}\quad Theorem \ref{Ai} ensures that
$$p_{d+k,d} = \sum_{\lambda \setminus \mu} (d+k)! \det \left( \frac{1}{(\lambda_i - \mu_j -i +j)!}\right)
$$
where the sum is over all skew shapes $\lambda \setminus \mu$ of
size $d+k$ having $k$ rows and such that two consecutive rows have
exactly two columns in common.

For such a skew shape $\lambda \setminus \mu$, let us define the
sequence $a = (a_1, a_2, \ldots , a_k)$ by $a_i = \lambda_i -
\mu_i$, $1\leq i \leq k$. The sequence $a$ is such that $a_1 ,a_2
,\ldots ,a_k \geq 2$ and $a_1 +a_2 +\cdots+a_k =d+k$. From the
remark at the end of the previous section, we additionally have
that $\mu_k = 0$ (point $4.$) and that $\mu_i = \lambda_{i+1} -2$
for all $1 \leq i \leq k-1$ (point $5.$). It is now trivial matter
to check that the sequence $a = (a_1, a_2, \ldots , a_k)$
completely and uniquely determines $\lambda \setminus \mu$. Hence,
the sum in the above formula can be taken over sequences $a =
(a_1, a_2, \ldots , a_k)$ such that $a_1 ,a_2 ,\ldots ,a_k \geq 2$
and $a_1 +a_2 +\cdots+a_k =d+k$.

In what follows, we give expressions of the entries
$a_{i,j}=\frac{1}{(\lambda_i -\mu_j -i+j)!}$ in terms of $(a_1,
a_2, \ldots , a_k)$.

If $i\leq j$,
then, by Theorem \ref{Ai}, $a_{i,j}=\frac{1}{(\lambda_i -\mu_j
-i+j)!}$. Thanks to the remarks stated at the end of the previous
section, we have that:
\begin{displaymath}
a_i +a_{i+1}+\ldots +a_j =(\lambda_i -\mu_i)+(\lambda_{i+1}
-\mu_{i+1})+\ldots (\lambda_j -\mu_j)=\lambda_i -\mu_j +2(j-i).
\end{displaymath}

This yields for the denominator of the above fraction the
following expression (leaving aside the factorial):
\begin{displaymath}
\lambda_i -\mu_j -i+j=a_i +a_{i+1}+\ldots +a_j +i-j,
\end{displaymath}
as desired.

If $j=i-1$, then we have immediately:
\begin{displaymath}
a_{i,i-1}=\frac{1}{(\lambda_i -\mu_{i-1} -1)!}=\frac{1}{(2-1)!}=1.
\end{displaymath}

Concerning the case $j=i-2$, since $\mu_{i-1} = \lambda_i -2$ and
$\mu_{i-2} \geq \mu_{i-1}$, we observe that $\lambda_i
-\mu_{i-2}\leq 2$, and that the equality holds precisely when
$\mu_{i-1} = \mu_{i-2}$, \emph{i.e.} when $a_i
=\lambda_{i-1}-\mu_{i-1}=2$. Thus we get:
\begin{displaymath}
a_{i,i-2}=\frac{1}{(\lambda_i
-\mu_{i-2}-i+i-2)!}=\chi_{a_{i-1}=2}.
\end{displaymath}

Finally, if $i>j+2$, then the denominator of $a_{i,j}$ is easily
seen to be negative, hence $a_{i,j}=0$.\cvd

From a theoretical point of view, Theorem \ref{principale}
completely solves the problem of the enumeration of $d$-minimal
permutations with respect to their length, giving a formula for
$p_{d+k,d}$. Unfortunately, it is clear that such a formula is
very difficult to use in concrete cases, due to its intrinsic
complexity. However, using our result we are able to rediscover
some known cases (namely $k=1,2$) and to give an interpretation of
Catalan numbers (corresponding to the case $k=d$), as well as to
get a formula for the case $k=3$ (that is, $d$-minimal
permutations of length $d+3$), which was first discovered in
\cite{MY} with different methods in terms of generating functions.
The formula we derive has been found with the help of Maple.

We start by collecting in a single theorem the known cases
$d=1,2$, showing how they can be derived from Theorem
\ref{principale}.

\begin{teor}(\cite{BP}) The following equalities hold:
\begin{eqnarray*}
p_{d+1,d}&=&1,
\\ p_{d+2,d}&=&2^{d+2}-(d+1)(d+2)-2.
\end{eqnarray*}
\end{teor}

\emph{Proof.}\quad When $k=1$, the formula of Theorem
\ref{principale} becomes completely trivial:
\begin{displaymath}
p_{d+1,d}=(d+1)!\cdot \left| \frac{1}{(d+1)!}\right|=1.
\end{displaymath}

In the case $k=2$, we have a single sum where a $2\times 2$
determinant appears:
\begin{eqnarray*}
p_{d+2,d}&=&\sum_{a_1 ,a_2 \geq 2\atop a_1 +a_2 =d+2}(d+2)!\cdot
\left| \begin{array}{cc} \frac{1}{a_1 !} & \frac{1}{(d+1)!}
\\ 1 & \frac{1}{a_2 !}
\end{array} \right| =\sum_{a=2}^{d}\left( {d+2\choose
a}-(d+2)\right)
\\ &=&2^{d+2}-2(d+3)-(d-1)(d+2)=2^{d+2}-2-(d+1)(d+2).
\end{eqnarray*}
\cvd

If $k=d$, the formula of Theorem \ref{principale} gives an
evaluation of Catalan numbers $(C_n )_{\ninN}$. The fact that
$p_{2d,d}$ is the $d$-th Catalan number is clear from Proposition
\ref{bij}, since $d$-minimal permutations are in bijection with
Young tableaux of rectangular shape having $d$ rows and 2 columns
(see \cite{St2}). Another
combinatorial proof of this fact is given in \cite{BP}. Thus we
get the following expression for Catalan numbers, which we have
not been able to find in the literature:
\begin{displaymath}
C_d =p_{2d,d}=(2d)!\cdot \left|
\begin{array}{ccccccc} \frac{1}{2!} & \frac{1}{3!} & \frac{1}{4!}
& \frac{1}{5!} & \frac{1}{6!} & \cdots & \frac{1}{(d+1)!}
\\ 1 & \frac{1}{2!} & \frac{1}{3!} & \frac{1}{4!} & \frac{1}{5!} & \cdots & \frac{1}{d!}
\\ 1 & 1 & \frac{1}{2!} & \frac{1}{3!} & \frac{1}{4!} & \cdots & \frac{1}{(d-1)!}
\\ 0 & 1 & 1 & \frac{1}{2!} & \frac{1}{3!} & \cdots & \frac{1}{(d-2)!}
\\ 0 & 0 & 1 & 1 & \frac{1}{2!} & \cdots & \frac{1}{(d-3)!}
\\ \vdots & \vdots & \vdots & \vdots & \vdots & \ddots & \vdots
\\ 0 & 0 & 0 & 0 & 0 & \cdots & \frac{1}{2!}
\end{array} \right| .
\end{displaymath}

We close the section with the evaluation of $p_{d+3,d}$. As we
stated above, to compute this value we have made extensive use of
Maple.

\begin{teor} The following equality holds:
\begin{displaymath}
p_{d+3,d}=3^{d+3}-(d^2 +4d+7)\cdot 2^{d+2} +\frac{1}{2}d^4
+\frac{5}{2}d^3 +\frac{33}{4}d^2 +6d-8.
\end{displaymath}
\label{d+3}
\end{teor}

\emph{Proof.}\quad We just have to apply Theorem \ref{principale}
in the case $k=3$, thus obtaining:
\begin{equation}\label{caso3}
p_{d+3,d}=\sum_{a,b,c\geq 2\atop a+b+c=d+3}(d+3)!\cdot \left|
\begin{array}{ccc}
\frac{1}{a!} & \frac{1}{(a+b-1)!} & \frac{1}{(d+1)!}
\\ 1 & \frac{1}{b!} & \frac{1}{(b+c-1)!}
\\ \chi_{b=2} & 1 & \frac{1}{c!}
\end{array} \right| .
\end{equation}

The presence of the characteristic function $\chi_{b=2}$ suggests
to consider two distinct cases.

\begin{itemize}

\item[i)] In Formula (\ref{caso3}), the partial sum for the tuples $(a,b,c)$ such that $b=2$ is:
\begin{eqnarray}\label{beq2}
&&\sum_{a,c\geq 2\atop a+c=d+1}(d+3)!\cdot \left|
\begin{array}{ccc}
\frac{1}{a!} & \frac{1}{(a+1)!} & \frac{1}{(d+1)!}
\\ 1 & \frac{1}{2} & \frac{1}{(c+1)!}
\\ 1 & 1 & \frac{1}{c!}
\end{array} \right| \nonumber
\\ &=&\sum_{a=2}^{d-1}(d+3)!\cdot \left|
\begin{array}{ccc}
\frac{1}{a!} & \frac{1}{(a+1)!} & \frac{1}{(d+1)!}
\\ 1 & \frac{1}{2} & \frac{1}{(d+2-a)!}
\\ 1 & 1 & \frac{1}{(d+1-a)!}
\end{array} \right| \nonumber
\\ &=&\sum_{a=2}^{d-1}\left( {d+3\choose a,2,d+1-a}+{d+3\choose
a+1}+(d+2)(d+3)+ \right. \nonumber
\\ &&\left. -\frac{1}{2}(d+2)(d+3)-(d+3){d+2\choose a}-(d+3){d+2\choose
a+1} \right) \nonumber
\\ &=&\sum_{a=2}^{d-1}\left( \frac{(d+2)(d+3)}{2}{d+1\choose a}+\frac{1}{2}(d+2)(d+3)-(d+2){d+3\choose
a+1} \right) \nonumber
\\ &=&{d+3\choose 2}(2^{d+1}-d-6)-(d+2)(2^{d+3}-d^2-7d-14).
\end{eqnarray}

\item[ii)] The partial sum for the tuples $(a,b,c)$ such that $b\neq 2$ is a bit more complicated to compute but gives the following:
\begin{eqnarray}\label{uffa}
&&\sum_{a,c\geq 2,b>2\atop a+b+c=d+3}(d+3)!\cdot \left|
\begin{array}{ccc}
\frac{1}{a!} & \frac{1}{(a+b-1)!} & \frac{1}{(d+1)!}
\\ 1 & \frac{1}{b!} & \frac{1}{(b+c-1)!}
\\ 0 & 1 & \frac{1}{c!}
\end{array} \right| \nonumber
\\ &=&\sum_{a,c\geq 2,b>2\atop a+b+c=d+3}\left( {d+3\choose
a,b,c}+(d+2)(d+3)-(d+3){d+2\choose a}-(d+3){d+2\choose c}\right)
\nonumber
\\ &=&\sum_{a=2}^{d-2}\sum_{c=2}^{d-a}{d+3\choose
a,d+3-a-c,c}+\sum_{a=2}^{d-2}\sum_{c=2}^{d-a}(d+2)(d+3)\nonumber
\\ &&-\sum_{a=2}^{d-2}\sum_{c=2}^{d-a}(d+3){d+2\choose a}-\sum_{a=2}^{d-2}\sum_{c=2}^{d-a}(d+3){d+2\choose
c}\nonumber
\\ &=&\alpha +\beta -2\gamma,
\end{eqnarray}
if we set
$$\alpha = \sum_{a=2}^{d-2}\sum_{c=2}^{d-a}{d+3\choose a,d+3-a-c,c}, \quad \beta = \sum_{a=2}^{d-2}\sum_{c=2}^{d-a}(d+2)(d+3)$$
$$\textrm{ and } \gamma = \sum_{a=2}^{d-2}\sum_{c=2}^{d-a}(d+3){d+2\choose a}.$$

We compute the two terms $\alpha$ and $\gamma$ using Maple:
\begin{eqnarray*}
\alpha &=& 3^{d+3}-(d+11)(d+6)\cdot 2^d +(d^3 +10d^2 +37d +51)
\\ \gamma &=& (d^2 -d-12)\cdot 2^{d+1}-\left( \frac{1}{2}d^3 -d^2
-\frac{41}{2}d-39\right) .
\end{eqnarray*}

Instead, the term $\beta$ is of course very easy to compute
directly:
\begin{displaymath}
\beta =(d+2)(d+3)\frac{(d-3)(d-2)}{2}=\frac{(d^2 -4)(d^2 -9)}{2}.
\end{displaymath}

Thus using Formula (\ref{uffa}) we get:
\begin{eqnarray}\label{bneq2}
&&\sum_{a,c\geq 2,b>2\atop a+b+c=d+3}(d+3)!\cdot \left|
\begin{array}{ccc}
\frac{1}{a!} & \frac{1}{(a+b-1)!} & \frac{1}{(d+1)!}
\\ 1 & \frac{1}{b!} & \frac{1}{(b+c-1)!}
\\ 0 & 1 & \frac{1}{c!}
\end{array} \right| \nonumber
\\ &=& 3^{d+3}-(5d^2 +13d+18)\cdot 2^d +\left( \frac{1}{2}d^4 +2d^3
+\frac{3}{2}d^2 -4d-9\right) .
\end{eqnarray}

\end{itemize}

Now, to finish the proof of our theorem, we just have to sum the
result of (\ref{beq2}) and (\ref{bneq2}), thus obtaining:
\begin{displaymath}
p_{d+3,d}=3^{d+3}-(d^2 +4d+7)\cdot 2^{d+2} +\frac{1}{2}d^4
+\frac{5}{2}d^3 +5d^2 +6d+1.
\end{displaymath}
\cvd

Table \ref{table:p_d+3,3} shows the first few terms of the sequence $(p_{d+3,3})_d$.

\begin{center}
\begin{table}[ht]
\begin{center}
\begin{tabular}{|r|c|c|c|c|c|c|c|}
\hline
$d$ & $3$ & $4$ & $5$ & $6$ & $7$ & $8$ & $9$ \\
\hline
$p_{d+3,3}$ & $5$ & $84$ & $686$ & $3936$ & $18387$ & $75372$ & $283052$\\
\hline
\end{tabular}
\end{center}
\caption{The first few terms of the sequence $(p_{d+3,3})_d$ \label{table:p_d+3,3}}
\end{table}
\end{center}

\section{A generalization} \label{section_generalization}

The main motivation of the present paper is the study of
$d$-minimal permutations and, in particular, their enumeration.
Our approach is based on a bijection between the set of $d$-minimal
permutations and a special class of skew Young tableaux, namely
those in which every pair of consecutive rows has precisely two
columns in common. Denote by $SkYT_2 (n,k)$ this set of tableaux,
$n$ being the number of cells and $k$ the number of rows. We can
generalize this setting in a very natural way, by defining the set
$SkYT_h (n,k)$ of skew Young tableaux having $n$ cells and $k$
rows such that any two consecutive rows have precisely $h\geq 1$
columns in common. In this final section we wish to relate these
tableaux with some families of permutations, as well as to describe
some enumerative results for low values of $h$.

\bigskip

Our first result is a generalization of Theorem \ref{bij}.

\begin{teor}\label{h-version} The set $SkYT_h (d+k,k)$ is in bijection with the set
of permutations of length $d+k$ having exactly $d$ descents and
satisfying the following property (call it $DES_h$):
$$DES_h
\begin{cases}
\text{for every $i\leq h-1$, if one deletes $i$ elements of a
permutation} \\ \text{and renames the remaining elements in the usual way,} \\
\text{the resulting permutation has precisely $d-i$ descents.}
\end{cases}
$$
\end{teor}

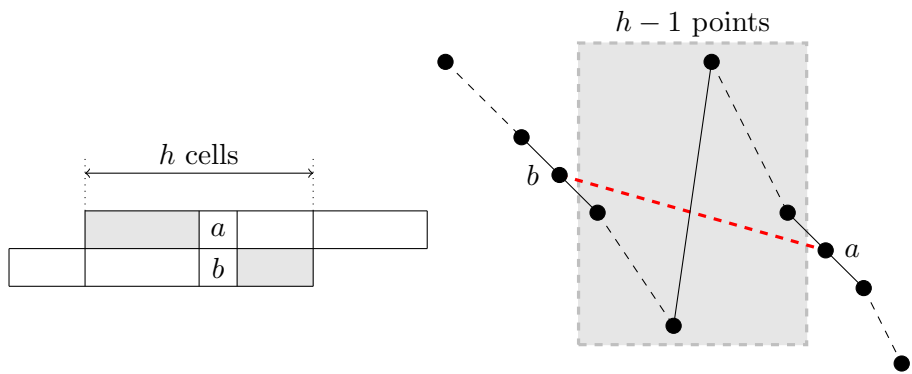
\begin{figure}[ht]
\begin{center}
\begin{tikzpicture}
\begin{scope}[scale=0.5]
\draw[fill,color=gray!20] (3,1) rectangle (6,2);
\draw[fill,color=gray!20] (7,0) rectangle (9,1);
\draw (6.5,0.5) node {$b$};
\draw (6.5,1.5) node {$a$};
\draw (1,0) -- (9,0);
\draw (1,1) -- (12,1);
\draw (3,2) -- (12,2);
\draw (1,0) -- (1,1);
\draw (3,0) -- (3,2);
\draw (6,0) -- (6,2);
\draw (7,0) -- (7,2);
\draw (9,0) -- (9,2);
\draw (12,1) -- (12,2);
\draw[dotted] (3,2) -- (3,3.5);
\draw[dotted] (9,2) -- (9,3.5);
\draw[<->] (3,3) -- (9,3);
\draw (6,3.5) node {$h$ cells};
\draw (2,-2) node {~};
\end{scope}
\end{tikzpicture} \begin{tikzpicture}
\begin{scope}[scale=0.5]
\draw[fill,color=gray!20] (4.5,1.5) rectangle (10.5,9.5);
\draw[color=gray!50, dashed, very thick] (4.5,1.5) rectangle (10.5,9.5);
\draw (7.5,10) node {$h-1$ points};
\draw[color=red, dashed, very thick] (4,6) -- (11,4);
\draw (1,9) [fill] circle (.2);
\draw[dashed] (1,9) --  (3,7);
\draw (3,7) [fill] circle (.2);
\draw (3,7) --  (4,6);
\draw (4,6) [fill] circle (.2);
\draw (3.3,6) node {$b$};
\draw  (4,6) -- (5,5);
\draw (5,5) [fill] circle (.2);
\draw[dashed] (5,5) -- (7,2) ;
\draw (7,2) [fill] circle (.2);
\draw (7,2) -- (8,9);
\draw (8,9) [fill] circle (.2);
\draw[dashed] (8,9) -- (10,5)  ;
\draw (10,5) [fill] circle (.2);
\draw (10,5) -- (11,4);
\draw (11,4) [fill] circle (.2);
\draw (11.7,4) node {$a$};
\draw (11,4) -- (12,3)  ;
\draw (12,3) [fill] circle (.2);
\draw[dashed] (12,3)  -- (13,1);
\draw (13,1) [fill] circle (.2);
\end{scope}
\end{tikzpicture}
\end{center}
\caption{Proof of Theorem \ref{h-version}. The gray area corresponds to the $h-1$ points removed in the proof. \label{fig:remove}}
\end{figure}

\emph{Proof.}\quad
Denote with $S_{d+k}^{(h)}(d)$ the set of permutations of length
$d+k$ having precisely $d$ descents and satisfying $DES_h$. Define
a map $f:S_{d+k}^{(h)}(d)\rightarrow SkYT_h (d+k,k)$ by suitably
generalizing the one given in Proposition \ref{bij}: starting from
the bottom, the $i$-th row of the tableau consists of the elements
of the $i$-th descending run of $\pi$, and two consecutive rows
are required to have exactly $h$ columns in common. We claim that
this map is well-defined. Indeed suppose, ab absurdo, that in
$f(\pi )$ there is a column in which $a$ is above $b$ and $a<b$. Without loss of generality, we can assume that $a$ is in the cell immediately above $b$.
Then, removing the $h-1$ entries of $\pi$ preceding $a$ and
following $b$ and belonging to the columns common to the rows of $a$ and
$b$, we obtain a permutation not satisfying $DES_h$ (the number of
descents is easily seen to be $d-h+2$), as shown on Figure \ref{fig:remove}. The fact that $f$ is
injective and surjective is trivial, and follows directly from its
definition.\cvd

%
%

It is clear that, when $h=2$, we get precisely Theorem \ref{bij},
since the resulting class of permutations is that of $d$-minimal
permutations.

\bigskip

We also have a characterization of the above classes of permutations in terms of patterns, which follows quite easily from the
above theorem, and so will be stated without proof.

\begin{teor} A permutation $\sigma$ belongs to $S_{d+k}^{(h)}(d)$,
for some $k$, if and only if it has exactly $d$ descents and its
ascents occur in the middle of a consecutive pattern of the form
$\pi =\pi_1 \pi_2$, where $ \pi_1$ and $\pi_2$ are words of the
same length $h$, both decreasing and $\pi_1 <\pi_2$
componentwise.
\end{teor}

Having introduced this generalized setting, it is natural to ask
what happens when $h<2$.

If $h=1$, what we obtain is the class of permutations having
exactly $d$ descents. It is well known that the number of
permutations of length $n$ having $d$ descents is given by the
Eulerian number $E_{n,d}$ (sequence A008292 in \cite{Sl}). Thanks
to our approach, we find a determinant expression of Eulerian
numbers which is believed to be new. Once again, the key
ingredient to obtain such a formula is of course Theorem \ref{Ai}.

\begin{teor}\label{secondario}  The number $E_{d+k,d}$ of
permutations of length $d+k$ having exactly $d$ descents (i.e.
satisfying condition $DES_1$) is
\begin{displaymath}
E_{d+k,d}=\sum_{a_1 ,a_2 ,\ldots ,a_k \geq 1\atop a_1 +a_2 +\cdots
+a_k =d+k}(d+k)!\cdot \det (B(a_1 ,\ldots ,a_k )),
\end{displaymath}
where $B(a_1 ,\ldots ,a_k )$ is the following matrix:
\begin{displaymath}
\left(
\begin{array}{cccccc}
\frac{1}{a_1 !}&\frac{1}{(a_1 +a_2 )!}&\frac{1}{(a_1 +a_2 +a_3
)!}&\cdots&\frac{1}{(a_1 +\cdots +a_{k-1})!}&\frac{1}{(a_1 +\cdots
+a_{k})!}
\\ 1&\frac{1}{a_2 !}&\frac{1}{(a_2 +a_3 )!}&\cdots&\frac{1}{(a_2 +\cdots +a_{k-1})!}&\frac{1}{(a_2 +\cdots +a_{k})!}
\\ 0&1&\frac{1}{a_3 !}&\cdots&\frac{1}{(a_3 +\cdots +a_{k-1})!}&\frac{1}{(a_3 +\cdots +a_{k})!}
\\ \vdots&\vdots&1&\ddots&\vdots&\vdots
\\ 0&0&\cdots&\ddots&\frac{1}{a_{k-1}!}&\frac{1}{(a_{k-1}+a_k )!}
\\ 0&0&0&\cdots&1&\frac{1}{a_k !}
\end{array}\right) .
\end{displaymath}

In other words, $B(a_1 ,\ldots ,a_k )$ is the $k\times k$ matrix
whose entries $b_{i,j}$ obey the following equalities:
\begin{eqnarray*}
a_{i,j}&=&\frac{1}{(a_i +\cdots +a_j )!},\qquad \textnormal{when
$i\leq j$,}
\\ a_{i,i-1}&=&1,
\\ a_{i,j}&=&0,\qquad \textnormal{when $i\geq j+2$.}
\end{eqnarray*}
\end{teor}

\emph{Proof.}\quad The proof essentially follows the same lines of
the proof of Theorem \ref{principale}; just observe that, in
this case, it is $\lambda_{i+1}-\mu_i =1$.\cvd

Moreover, the determinant of the matrix $B(a_1 ,\ldots ,a_k )$ has a very neat
recursive expression, from which a closed formula can be deduced.

\begin{prop} Set $D(a_1 ,\ldots ,a_k )=\det (B(a_1 ,\ldots ,a_k
))$. Then
\begin{displaymath}
D(a_1 ,a_2 ,\ldots ,a_k )=\frac{1}{a_1 !}\cdot D(a_2 ,a_3 ,\ldots
,a_k )-D(a_1 +a_2 ,a_3 ,\ldots ,a_k ).
\end{displaymath}
\end{prop}

\emph{Proof.}\quad Just expand $D(a_1 ,\ldots ,a_k )$ with respect
to its first column.\cvd

\begin{cor} The following formula holds:
\begin{equation}\label{determinant}
D(a_1 ,\ldots ,a_k )=\sum_{i=1}^{k}(-1)^{k-i}\cdot \sum_{\alpha
=(\alpha_1 ,\ldots ,\alpha_i )\atop \alpha \in
\mathbf{PL}(a_1,\cdots ,a_k)}\frac{1}{|\alpha_1 |!\cdot \ldots
\cdot |\alpha_i |!},
\end{equation}
where $\mathbf{PL}(a_1,\cdots ,a_k)$ denotes the set of linear
partitions of the totally ordered set $\{ a_1 ,\ldots ,a_k \}$ and
$|\alpha_i |$ is the sum of the elements of the block $\alpha_i$.
\end{cor}

\emph{Proof.}\quad We start by observing that, when $k=1$, the
outer sum of the r.h.s of (\ref{determinant}) reduces to a single
summand (for $i=1$), as well as the inner sum, which has the
unique summand $\frac{1}{a_1 !}$. Moreover, when $k=2$, the r.h.s.
of (\ref{determinant}) consists of two summands, which are
$(-1)\cdot \frac{1}{(a_1 +a_2 )!}$ (for $i=1$) and $\frac{1}{a_1
!}\cdot \frac{1}{a_2 !}$ (for $i=2$), and this coincides with the
expression of $D(a_1 ,a_2 )$.

We can now conclude our proof using an inductive argument. The set
$\mathbf{PL}(a_1,\cdots ,a_k)$ can be partitioned into two
subsets, namely the linear partitions in which $a_1$ occurs as a
singleton (call this subset $X$) and the linear partitions in
which $a_1$ occurs in a block of cardinality at least 2 (call this
subset $Y$). Using this partition of $\mathbf{PL}(a_1,\cdots
,a_k)$ we can split the sum in the r.h.s of (\ref{determinant})
into two sums, the first taking into account the contribution of
$X$ and the second taking into account the contribution of $Y$. We
thus obtain the following equalities:
\begin{eqnarray*}
\sum_{\alpha =(\alpha_1 ,\ldots ,\alpha_i )\atop \alpha \in
\mathbf{PL}(a_1,\cdots ,a_k)}\frac{1}{|\alpha_1 |!\cdot \ldots
\cdot |\alpha_i |!}&=&\frac{1}{a_1 !}\cdot \sum_{\beta =(\beta_1
,\ldots ,\beta_j )\atop \beta \in \mathbf{PL}(a_2,\cdots
,a_k)}\frac{1}{|\beta_1 |!\cdot \ldots \cdot |\beta_j |!}\\
&+&\sum_{\gamma =(\gamma_1 ,\ldots ,\gamma_t )\atop \gamma \in
\mathbf{PL}(\{ a_1 ,a_2 \},a_3 ,\cdots ,a_k)}\frac{1}{|\gamma_1
|!\cdot \ldots \cdot |\gamma_t |!},
\end{eqnarray*}
whence, using the induction hypothesis and the above proposition:
\begin{eqnarray*}
\sum_{i=1}^{k}(-1)^{k-i}\cdot \sum_{\alpha =(\alpha_1 ,\ldots
,\alpha_i )\atop \alpha \in \mathbf{PL}(a_1,\cdots
,a_k)}\frac{1}{|\alpha_1 |!\cdot \ldots \cdot |\alpha_i
|!} \qquad \qquad \qquad \qquad \qquad \qquad \quad \\ = \frac{1}{a_1 !}\cdot D(a_2 ,\ldots ,a_k)-D(a_1 +a_2 ,a_3
,\ldots ,a_k ) = D(a_1 ,\ldots ,a_k),
\end{eqnarray*}
as desired.\cvd


\emph{Remark.}\quad An alternative approach to the case $h=1$
could be done via the notion of Hessenberg matrix. An (upper)
\emph{Hessenberg matrix} is a square matrix having zero entries
below the first subdiagonal. Hessenberg matrices prove their
usefulness especially in numerical analysis and computer
programming, being a sort of normal form to which any square
matrix can be reduced in a finite number of steps. There are also
some papers in the literature concerning the evaluation of the
determinant of certain Hessenberg matrices having special form
(see for instance \cite{BS} and \cite{LCT}). In \cite{T}, the
determinant of Hessenberg matrices having all the elements of the
first subdiagonal equal to 1 is considered (this is precisely the
kind of matrices we meet in Theorem \ref{secondario}).

\bigskip

Theorem \ref{h-version} does not have meaning when $h=0$. The
corresponding set $SkYT_0 (n,k)$ consists of all skew Young
tableaux having $n$ cells and $k$ rows such that any two
consecutive rows only have the corners of two cells in common. In
this case, it is immediate to see that $SkYT_0 (n,k)$ is in
bijection with all surjective functions from an $n$-set to a
$k$-set: just interpret the elements of a tableau as balls and the
rows of a tableau as boxes. Thus we get immediately that $|SkYT_0
(n,k)|=k!\cdot S(n,k)$, where the $S(n,k)$'s are the Stirling
numbers of the second kind.
We can also use Theorem \ref{Ai} to get an
analog of Theorems \ref{principale} and \ref{secondario}; indeed,
we can derive the following formula:
\begin{displaymath}
|SkYT_0 (d+k,k)|=\sum_{a_1 ,a_2 ,\ldots ,a_k \geq 1\atop a_1 +a_2
+\cdots +a_k =d+k}(d+k)!\cdot \det (C(a_1 ,\ldots ,a_k )),
\end{displaymath}
where $C(a_1 ,\ldots ,a_k )$ is the following triangular matrix:
\begin{displaymath}
\left(
\begin{array}{cccccc}
\frac{1}{a_1 !}&\frac{1}{(a_1 +a_2 +1)!}&\frac{1}{(a_1 +a_2 +a_3 +2
)!}&\cdots&\frac{1}{(a_1 +\cdots +a_{k-1} +k-2)!}&\frac{1}{(a_1 +\cdots
+a_{k} +k-1)!}
\\ 0&\frac{1}{a_2 !}&\frac{1}{(a_2 +a_3 +1)!}&\cdots&\frac{1}{(a_2 +\cdots +a_{k-1} +k-3)!}&\frac{1}{(a_2 +\cdots +a_{k}+k-2)!}
\\ 0&0&\frac{1}{a_3 !}&\cdots&\frac{1}{(a_3 +\cdots +a_{k-1}+k-4)!}&\frac{1}{(a_3 +\cdots +a_{k}+k-3)!}
\\ \vdots&\vdots&\ddots&\ddots&\vdots&\vdots
\\ 0&0&\cdots&0&\frac{1}{a_{k-1}!}&\frac{1}{(a_{k-1}+a_k +1)!}
\\ 0&0&0&\cdots&0&\frac{1}{a_k !}
\end{array}\right) .
\end{displaymath}

From here it is then immediate to obtain $|SkYT_0
(d+k,k)|=\sum_{a_1 ,a_2 ,\ldots ,a_k \geq 1\atop a_1 +a_2 +\cdots
+a_k =d+k}{d+k\choose a_1 ,\ldots a_k}$ which is known to be the
number of surjective functions from an $(d+k)$-set to a $k$-set,
as already shown a few lines above.

\section{Further work}

Even if our approach to the enumeration of $d$-minimal
permutations allows us to completely solve the problem from a
purely theoretical point of view, it is doubtless that its
application to concrete cases shows some technical difficulties.
This is of course due to the intrinsic complexity of the sums of
determinants appearing in Theorem \ref{principale}. However, it
seems plausible that at least a few more cases than those we deal
with in the present paper can be managed by means of our
technique. 

Another interesting problem that remains untouched concerns the
study of the structure of the poset determined by a minimal
permutations with $d$ descents, defined in \cite{BP} and recalled
in Section \ref{prel} here. For instance, one can observe that a
$d$-minimal permutation corresponds to a linear extension of the
associated poset. Moreover, an interesting (and classical) line of
research could be the investigation of the properties of the
distributive lattice of the sup-irreducibles of these posets.

\end{document}